\documentclass{article}
\usepackage[utf8]{inputenc}

\begin{document}

\title{Maximum Relative Divergence Principle for Grading Functions on Direct Products of Chains}

\author{Alexander Dukhovny \\
\small{Department of Mathematics, San Francisco State University} \\
\small{San Francisco, CA 94132, USA} \\
\small{dukhovny@sfsu.edu} \\ }

\date{\today}          

\maketitle

\begin{abstract}

The concept of Shannon Entropy for probability distributions and associated Maximum Entropy Principle are extended here to the concepts of Relative Divergence of one Grading Function from another and Maximum Relative Divergence Principle for grading functions on direct products of totally ordered chains (chain bundles). Several Operations Research applications are analyzed.

\end{abstract}

\section{Introduction}
\label{sec:Intro}

In numerous probability theory problems and beyond The Insufficient Reason Principle has been used in the form Maximum (Shannon) Entropy Principle (MEP). According to MEP, the "most reasonable" (using the fewest extra assumptions) way to determine missing pieces of the needed probability distribution is to maximize, under some application-specific constraints, the Shannon Entropy functional (see, e.g., \cite{Shannon}). That approach has proved effective in so many cases that references are just too many to quote.\par

The original Shannon Entropy formula and MEP have been generalized in many contexts. The list includes Relative Entropy, Kullback-Leibler Divergence, Partition Entropy,  Kolmogorov-Sinai Entropy, Topological Entropy, Entropy of general non-probabilistic measures (capacities) and great many others  (see, e.g., references [5-11] and a review in \cite{Review}.

In our initial paper \cite{AOME}, Shannon Entropy was generalized to the concept of Relative Divergence (RD) of one Grading Function (GF) from another on a totally ordered set (chain), reducing to Shannon Entropy in a special case. (The term Relative Divergence was chosen in keeping with Kullback-Leibler Divergence - see \cite{K-L} - also known as Relative Divergence of probability measures.)

Further, in \cite{MRDP on P-sets} we started the process of extending the concepts of Relative Divergence and Maximum Relative Divergence Principle (MRDP) to partially ordered sets. In the process, we demonstrated that 

1. both of those concepts reduce to Shannon Entropy and Maximum Entropy Principle when used to determine the "most reasonable" probability distribution on a sample space.

2. conclusions made by using MRDP in new problems agree with "common sense" ones - where such are available;

3. MRDP can be effectively used in new applications.

The new applications in \cite{MRDP on P-sets} included working with a power set $W = 2^X$ of the event space $X$ ordered by subset inclusion. In that case our results coincided with the ones obtained in our earlier work using the concept of General Entropy of General (non-additive) Measures (see \cite{DFuzzy}) for power sets. 

Namely, using a "normalized" General Measure $\mu (w), \quad w \in W$, as a Grading Function $F(w)$ of a subset $w$, it followed that its relative divergence from the "subset cardinality" grading function $N(w)$ reduces to the minimum of Shannon Entropy values of all probability distributions said to be "subordinate" to $\mu (w)$.

Applying MRDP in \cite{MRDP on P-sets} beyond Probability Theory applications made it possible to consider some problems arising in Operations Research. It showed that MRDP does lead to "most reasonable, common sense, natural" results where such results are available. Also, some new results were obtained under new application-specific constraints.

In this paper we continue the process started in \cite{AOME} and extend the RD  and MRDP concepts to sets that are direct products of totally ordered chains referred to as "chain bundles".

To make the paper self-contained, some basic general definitions and properties of RD and MRDP are presented in section \ref{Basics}, along with relevant technical tools from the Shannon Entropy theory. 

In section \ref{RD and MRDP on Bundles} general results for RD and MRDP are specified for direct products (bundles) of totally ordered chains. Based on single-chain cases, we explore possible solution process domain (the set of admissible grading functions) issues under external linear constraints.

Next, in sections \ref{Height-dependent GFs}, \ref{+ separable GFs} and \ref{MRDP for controlled GFs} we consider a number of special cases that emerge naturally in representative Operations Research applications. Using MRDP, we explore the MRDP solution process under several types of constraints on the admissible grading functions: imposed direct linear relations on their values and implied structural forms (referred to as "height-dependent", "additively separable", "parameters-controlled").

\section{Basic Definitions and properties}\label{Basics}

The initial setup of Relative Divergence in \cite{AOME} begins as follows: let $W$ be a set totally ordered by the order relation $\prec$. A real-valued function $F$ on $W$ is said to be a Grading Function (GF) on $W$ if it is order-monotonic, that is,

$w \prec v \iff F(w) < F(v)$ for all $w, v \in W $.

In this paper we consider a discrete countable 
 $W = \{ w_k, \quad k= \ldots, -1, 0, 1 , \ldots \}$
and refer to its "ordinal" function $I: I(w_k) = k$ as a "natural" GF on $W$ .

When $F(w)$ and $G(w)$ are GFs on a discrete countable $W$, the RD of $F$ from $G$ on $W$ is defined (assuming absolute convergence of the series), as

\begin{equation} \label{discrete RD}
\mathcal{D}(F \Vert G) \vert_W = -\ \sum_{k= - \infty}^\infty\
\ln  \left( \frac{\Delta_k F}{\Delta_k G}\right) \Delta_k F , 
\end{equation}

\noindent 
where \quad $\Delta_k F = F(w_k) - F(w_{k-1}), \quad 
\Delta_k G= G(w_k) - G(w_{k-1}), \quad  k = \ldots,-1,0,1,\ldots$.

It follows directly from the definition that $\mathcal{D}(F \Vert G)$ possesses some special properties w.r.t. linear transformations of those grading functions. 

1. $\mathcal{D}((c+F) \Vert (c+G)) \vert_W = \mathcal{D}(F \Vert G) \vert_W $.

2. $\mathcal{D}(cF \Vert cG) \vert_W =
c\mathcal{D}(F \Vert G) \vert_W, \quad \forall c>0.$

Also, when a GF $F$ is bounded on  $W$, denoting  by $m$ and $M$, respectively, the minimum and maximum of $F$ on $W$, it follows directly from (\ref{discrete RD}) that $\forall c > 0$

\begin{equation} \label{scaling F}
\mathcal{D}(cF \Vert I) \vert_W = 
 c\mathcal{D}(F \Vert I) \vert_W - c(M-m)\ln{c}, 
\end{equation}

\noindent where $\Delta_W F = M - m$.

Any bounded GF can be presented as $F(w) = m + \hat{F} (w) (M-m)$, where we define $\hat{F}:= \frac{F-m}{M-m}$.

Obviously, $\hat{F}$ is itself a grading function with the grading range [0,1] (the "standardized" $F$) which can be interpreted as a cumulative probability distribution function on $W$. Based on that, we will use the term  "Disorder Entropy" for $\hat{F}$ on $W$ and denote it $\mathcal{H}(\hat{F})\vert_{W}$:

\begin{equation} \label{H of F-hat}
\mathcal{H}(\hat{F})\vert_W := \mathcal{D}(\hat{F}\Vert I)\vert_W
\end{equation}

\noindent Now, using $c = M-m$ in (\ref{scaling F}) for $\hat{F}$, it follows that

\noindent $\mathcal{D}(F \Vert I) \vert_W = 
\mathcal{D}((M -m)\hat{F}) \Vert I) \vert_W =
(M-m) \mathcal{H}(\hat{F}) - (M-m) \ln{(M-m)}$,

\noindent which establishes the connection between the concepts of Relative Divergence and Shannon Entropy. As such, some classical techniques from the standard toolkit of Shannon Entropy theory (see, e.g., \cite{Jaynes}) become relevant here.

\textbf{Lemma 1.} For a probability distribution 

$\{p_i\}, \quad 0\leq p_i\leq 1, \quad \sum_{i=1}^n p_i = 1, \quad i = 1, \ldots, n$,

\noindent the maximum value of Shannon Entropy of that distribution

$\mathcal{H}\quad  = -\sum_{i=1}^n p_i\ln{p_i} \quad = \ln{n}$

\noindent is attained when $p_i = \frac{1}{n}, \quad i = 1, \ldots, n$. 

Based on Lemma 1, the following result from \cite{MRDP on P-sets} will be important in what follows.

\textbf{Lemma 2.} Let $F(i)$ be a grading function on the chain $W = \{0, \ldots, n \}$, and suppose some values of $F(i)$ are specified, that is,

$F(n_k)=M_k, \quad k=1, \ldots, K$, 

\noindent where we define 
$n_0 = 0, \quad M_0 = m,\quad n_K = n, \quad M_n = M$ 

\noindent and also \quad 
$\Delta_k M = M_k - M_{k-1}$, $\Delta_k n = n_k - n_{k-1}$, $k = 1, \ldots, K.$

Then the maximum value of 
$\mathcal{D}(F \Vert I) \vert_W = 
\sum_{k=1}^K [\Delta_k M \Delta_k n)- \Delta_k M \ln \Delta_k n)]$

\noindent is attained when $F$(i) is a piece-wise linear function:

\begin{equation} \label{Lemma 2 formula}
F(i)=a_k+b_k i, \quad i \in I_k, \quad k=1, \ldots, K, \quad \forall i \in W,
\end{equation}

\noindent where
$\quad b_k= \frac{ \Delta_k M} {\Delta_k n}, \quad a_k = M_k - b_k n_{k-1},$ 

\noindent and index intervals \quad 
$I_k = ( n_{k-1}, n_k ], \quad k = 1, \ldots, K$.

\textbf{Proof.}  Let

$q_{i,k} = \frac{F(i) - F(i-1)}{ M_k - M_{k-1}}, \quad i \in I_k$.

Maximizing $\mathcal{D}(F \Vert I) \vert_W$ in this case reduces to the following problem.

Find the $q_{i,k} \geq{0}, \quad i \in I_k, \quad k = 1, \ldots, K$

\noindent that maximize 
$ \quad -\sum_{k=1}^K \sum_{i \in I_k} q_{i,k} \ln (q_{i,k} ) $

\noindent subject to 
$\quad \quad \sum_{i \in I_k} q_{i,k} = 1, \quad k = 1, \ldots,K$.

The additive form of the maximized expression leads to $K$ independent maximization problems for each $k = 1, \ldots, K$. Using Lemma 1, the (unique) solution of each one of them is

$q_{i,k} = \frac{1}{\Delta_k n}, \quad i \in I_k$,

\noindent which completes the proof of Lemma 2.

In particular, when only $M_n =  M, \quad M_0 = m = 0$  are specified,  $K=1$, so (\ref{Lemma 2 formula}) reduces to

$F(i) = M \frac{i}{n},\quad \forall i \in W$,

$\mathcal{D}(F \Vert I) \vert_W = M \ln{n} - M \ln{M}$

\section{Relative Divergence and Maximum Relative Divergence Principle on Chain Bundles}\label{RD and MRDP on Bundles}

 In this section we specify general results for Relative Divergence of grading functions and Maximum Relative Divergence Principle to an "event space" $W$ - referred heretofore to as a chain bundle - a direct product of $R$ totally ordered chains: 
 
 $W = X_1 \times, \ldots, \times X_R$, where  $X_r = \{ x_r (i), \quad i = 0, 1, \ldots, n_r \}, \quad r = 1, \ldots, R$. 
 
Using vector notation, $W = \{ \vec{w} (\vec{i}) \} $
where the elements of $W$ are denoted as

$ \vec{w} (\vec{i}) = [x_1 (i_1), \ldots, x_N (i_R)]$

\noindent and their vector indices $\vec{i} = [i_1, \ldots, i_R]$.
 
As a direct product, the standard order relation of the elements of $W$ is imposed by the order of their vector indices: for all unequal vectors in $W$
 
 $\vec{w}_{\vec{i}} \prec  \vec{w}_{\vec{j}} \iff \vec{i} \prec \vec{j}$,
 
\noindent (that is, $i_t \leq j_t, \quad \forall t.$)

 As such, $W$ has the minimal element 
 $\vec{w}_{min} = [x_1 (0), \ldots, x_R (0)]$ 
 and the maximal element  
 $\vec{w}_{max} = [x_1 (n_1), \ldots, x_R (n_R)]$ 
 
 Two elements of $W$ are said to be adjacent if their vector indices differ by only one component where the difference is 1.

 A sequence of vectors of $W$ is said to be a maximal chain $MC$ in $W$ if it is totally ordered and no other chain in $W$ contains it. Consequently, all maximal chains in $W$ have the same minimal and maximal elements and the same number of elements $K = 1 + n_1 + \ldots + n_R$, and each element of $W$ belongs to at least one maximal chain. 
 
 Therefore, for any grading function $F(\vec{w})$ on $W$ its reduction to any maximal chain $MC \in $ has the same minimum and maximum values $m, M$ and overall grade spread $M-m$ on each $MC$.
 
 Because of that, following the approach of \cite{MRDP on P-sets}, and treating $W$ as a union of all of its maximal chains, we define the relative divergence of two grading functions on $W$ as follows:
 
 \begin{equation} \label{RD on W}
\mathcal{D}(F\Vert G)\vert_W := 
\min_{MC\subset W} \mathcal{D}(F\Vert G)\vert_{MC}
\end{equation}

(That definition follows the one proposed in \cite{DFuzzy} to facilitate Maximum Entropy Principle for a general (non additive) measure $\mu$ on the powerset $W$ of the element set $X$: its Shannon entropy $\mathcal{H} (\mu) \vert_W$ should be taken as the minimum over all maximal chains $MC$ in $W$ of Shannon Entropy values of probabilistic measures $\mu_{MC}^s$ said to be "subordinate to $\mu$ on $MC$". Each such $\mu_{MC}^s$ is completely and uniquely determined by its values on the subsets comprising $MC$ set equal to the values of $\mu$ on those subsets.)

Since the definition of RD does not involve actual values of elements of $W$, to simplify notation we will, where feasible, refer to the elements 
$\vec{w} (\vec{i)}$ simply by their index vectors $\vec{i}$.

Each maximal chain $MC$ in $W$ is a sequence of adjacent vectors 

$\{\vec{i} (k), \quad k=0, 1, \ldots, K \}$. 

For a grading function $F(\vec{i})$ defined on $W$ we also define 

$f_{MC} (\vec{i}(k)) = 
F(\vec{i} (k)) - F(\vec{i} (k-1)), \quad k = 1, \ldots, K,$ 

\noindent - the "increment" function of $F$ along the chain $MC$. As a grading function, $F$ is monotonic, so its increment function assumes only nonegative values.

In some applications it is possible to identify an obvious "common sense", "natural" GF on $W$. When $W$ is a direct product of chains, in the absence of constraints, the natural one is $N(\vec{i})$ - the "height" of $\vec{i}$ w.r.t. the order on $W$), the sum of all components of its vector index $\vec{i}$:

$N(\vec{i}) := i_1 + \ldots + i_R$.

Its increment function along any maximal chain $MC$ has the same constant values:

$f_{MC} (\vec{i}(k)) = 1, \quad k = 1, \ldots, K$.

Just as in our previous papers on the subject, we can now extend Maximum Entropy Principle (MEP) for probability distributions to Maximum Relative Divergence Principle (MRDP), as stated in \cite{AOME}, applied to grading functions on chain bundles (direct products of chains).

MRDP: An "admissible" (satisfying the constraints of the problem) bounded grading function $F$ on a chain bundle $W$ is said to be "the most reasonable" (within the fixed grading interval) if it gives maximum to 
$\mathcal{D}(F \Vert N) \vert_ W $.

In its most general setup, obtaining $\mathcal{D}(F\Vert G)\vert_W$ has a very high computational cost as the number of maximal chains in a bundle of $R$ chains increases rapidly with the number and the sizes of the bundled chains. 

In addition, using MRDP involves optimization analysis of all
$\mathcal{D}(F \Vert N) \vert_{MC}$  over all $MC$ in $W$,
possibly under application-specific constraints, increasing the overall computational load. In the next sections we explore a number of special cases that arise naturally in some applications under application-specific assumptions which reduce the computational cost of finding Relative Divergence and using MRDP.

At the same time, in some special cases, the nature of the problem of interest and the imposed constraints complicate the set of admissible grading functions on $W$ and the analysis of the MRDP problem.

For a simple illustration, consider the case where $W$ is a single chain: 
$W = \{ 0, 1, \ldots, n \} $, and an admissible grading function's $F(i)$ increments \quad $f_i = F(i) - F(i-1) \geq 0, \quad i = 1,\ldots, n,$ 
must satisfy two constraints: one imposed by the very definition of increments, and the other one that imposes a specified value of a linear combination of all $f_i$. Assuming for simplicity that
$F(0) = m = 0, \quad F(n) = M$, the constraints look as follows:

$\sum_{i=1}^n f_i = M$, \quad and \quad $\sum_{i=1}^n c_i f_i = \mu$.

Following MRDP in this case calls for maximizing the value of

$\mathcal{D}(F \Vert N) \vert_ W  = -\sum_{i=1}^n f_i \ln{f_i}$

\noindent under those constraints.

It is clear from the constraints that, depending on the coefficients of the linear constraint, the domain of the problem may be empty - as, for example, when either all $c_i > \frac{\mu}{M}$, or all $c_i < \frac{\mu}{M} $. 

When $c_i = \frac{\mu}{M}, \quad \forall i$, the constraints coincide and the MRDP problem reduces to the one covered by Lemma 2, so $F(i) = i\frac{M}{n}, \quad i = 0, 1, \ldots, n$.

Otherwise, when the maximum value of $-\sum_{i=1}^n f_i \ln{f_i}$ is attained at an interior point of the domain, using Lagrange multipliers, all $f_i$ of that solution, in addition to the constraints of the problem, must satisfy the following Lagrange equations:

$ -\ln{f_i} - 1 -\alpha - c_i \beta = 0, \quad i = 1, \ldots, n.$

Denoting $a = e^{-1-\alpha}$ and $b = e^{-\beta}$, it follows then that 

\begin{equation}\label{one chain, mu}
f_i = a b^{c_i}, \quad i = 1, \ldots,n,
\end{equation}

\noindent using which, the constraints can be restated as

$a\sum_{i=1}^n b^{c_i} = M$,

$a\sum_{i=1}^n c_ib^{c_i} = \mu$,

\noindent and used to specify $a$ and $b$.

Taking a representative case from the Probability Theory, consider a random variable assuming values in $W = \{ 1, \ldots, n \}$. Looking for its "most reasonable" probability distribution function $F(i)$ under the assumption that the expected value $\mu > 1$ is specified, one can treat $F(i)$ as a grading function on $W$, where $M=1, m=0$. The specified $\mu$ presents as a linear constraint where $c_i = i, \quad i=1, \ldots, n$. Using MRDP, it follows from (\ref{one chain, mu}) that the distribution must be geometric, exactly as it follows from MEP in the Probability Theory.

\section{"Height-dependent" Grading Functions}
\label{Height-dependent GFs}

When $W$ is a bundle of more than one chain $(R>1)$, using MRDP may still lead to smaller computational cost under some special assumptions on the application-imposed form of the grading functions, constraints and the structure of the involved chains.

When the application model imposes no constraints on the grading functions, then the MRDP-suggested most reasonable grading function must be linear: 

$F(\vec{i}) = m + N(\vec{i}) \frac{(M-m)}{K}$.

Indeed, in that case, by Lemma 2, for any maximal chain $MC$ in $W$, using that formula for $F$ results in the highest possible value of

$\mathcal{D}(F \Vert N) \vert_{MC} = (M-m)\ln K - (M-m)\ln{(M-m)}.$

At the same time, direct results can also be obtained under the constraint on the very nature of the admissible grading functions: $F(\vec{i})$ on $W$ must be "height-dependent", that is, 

$F(\vec{i}) = F(N(\vec{i}))$. 

Indeed, in that case, the formula of (\ref{RD on W})  yields

\begin{equation} \label{RD of F(N)}
\mathcal{D}(F \Vert N)\vert_W = - \sum_{k=1} ^{K} f(k) \ln {(f(k))},
\end{equation}

\noindent where $f(k) = F(k)-F(k-1), \quad k = 1, \ldots, K$.

A representative example arises in Queuing Theory: a service batch is formed by the server by selecting groups of customers from several waiting lines. The service batch is therefore an element of the bundle of the queues. The total batch service "cost" (say, total batch service time, including loading and processing), is often assumed to depend only on the total size of the batch. Assuming that the cost increases with the batch size, it is a "height-dependent" grading function on $W$.

By Lemma 2, in the absence of constraints, the MRDP-suggested most reasonable grading function in this case must be a linear function of the "height" of the element of $W$: 

$F(\vec{i}) = m + N(\vec{i}) \frac{(M-m)}{K}$.

When an admissible "height-dependent" grading function must, in addition, have its increments satisfy a linear constraint, the previous section single-chain analysis applies and yields direct results. 

\section{"Additively separable" Grading Functions}
\label{+ separable GFs}

In many applications where the event space $W$ of its model is a direct product (bundle) of chains, admissible grading functions are to be construed out of grading functions defined on the chains bundled in $W$.

Specifically, consider a direct product of two (or, consequently, any number) chain bundles (of different totally ordered chains):

\noindent $U = X_1 \times \ldots \times X_R$, 
\quad $V = Y_1 \times \ldots \times Y_Q$, and 

\noindent $W = U \times V = \{\vec{w} (\vec{s}) \} $, 
where vector index $\vec{s} = [\vec{i}, \vec{j}]$ 

\noindent and vector elements of $W$ are
$\vec{w} (\vec{s}) = [\vec{u} (\vec{i}), \vec{v} (\vec{j})]$.

A grading function $F_W$ on $W = U\times V$ is said to be "additively separable" on $W$ w.r.t. to $(U, V)$ if its values are sums of values of grading functions defined on $U$ and $V$, that is for any $\vec{w} = [\vec{u}, \vec{v}]$

\begin{equation}\label{add separable}
F_W(\vec{w}) = F_U(\vec{u}) + F_V(\vec{v})
\end{equation}

In particular, the "natural" GF $N$ on a chain bundle is additively separable:

$N_W(\vec{w}) = i_1 + \ldots + i_R + j_1 + \ldots j_Q = 
N_U(\vec{u}) + N_V(\vec{v)})$

An important representative example arises in the context of a Queuing Theory application: a server forms a service batch by seperately selecting groups of items from two separate waiting lines $U$ and $V$, in order of their lines, independently of the other line choice. The service batch is an element of $W=U \times V$. 

The total batch service "cost" (say, batch service time), in its simplest model is the sum of both groups' service costs, determined, in turn, by each group's size. Under a common assumption that those costs increase with the group size, the total service batch cost $F_W$ is a grading function on $W$ computed as a sum of the batch forming groups' costs $F_U$ and $F_V$. 

 Proposition 1. Suppose $F$ and $G$ are additively separable grading functions on chain bundles $U$ and $V$ and $W = U\times V$. Then

\begin{equation} \label{RD +separable}
\mathcal{D}(F \Vert G) \vert_W = 
\mathcal{D}(F_U \Vert G_U) \vert_U + \mathcal{D}(F_V \Vert G_V) \vert_V.
\end{equation}

Proof. The general definition formula in \ref{add separable} is a sum of terms involving involving only the increments of $F$ and $G$ corresponding to each of the values of $k, \quad k = 1,\ldots,K$ for each maximal chain $MC$.

Each $MC$ in $W = U \times V$ can be split into subchains $\{ MC_{\vec{i}(k)} \}$ in $V$ and $\{ MC_{\vec{j}(k)} \}$ in $U$, where $\vec{i}(k)$ and ${\vec{j}(k)}$, respectively, are constant. Clearly, both kinds of subchains spliced together form maximal chains $MC_V$ in $V$ and $MC_U$ in $U$, respectively.

Combining the terms of into groups corresponding to $MC_V$ in $V$ and $MC_U$ in $U$ and using the additive separability of both $F$ and $G$ one observes that in those groups the increments of those functions are either the increments of $F_V$ and $G_V$ or the increments of $F_U$ and $G_U$, respectively. As such, it follows that

$\mathcal{D}(F \Vert G) \vert_{MC} = 
\mathcal{D}(F_U \Vert G_U) \vert_{MC_U} + 
\mathcal{D}(F_V \Vert G_V) \vert_{MC_V}$.

Clearly, so obtained chains $MC_V$ in $V$ and $MC_U$ in $U$ contain all possible maximal chains in $U$ and $V$. As such, taking the minimum of the expression above over all $MC$ in $W$ completes the proof of Proposition 1.

Proposition 1 may simplify the way MRDP applies to the combined problem. Namely, if the constraints of the maximization problem are stated separately and independently for $F_U \vert_U$ and $F_V \vert_V$ the overall maximization problem splits into separate smaller maximization problems, resulting in massive reduction of the computing cost. Consequently, it opens a way to analyze an MRDP problem on the entire chain bundle by first identifying separate independent sub-bundles and splitting the combined problem into two (or more) smaller problems.

Corollary 1. If grading functions 
$F$ and $G$ on $W = X_1 \times \ldots \times X_R$ 
are completely additively separable, that is,

$F(\vec{w}) =  F_1(x_1) + \ldots + F_R(x_R)$,

\noindent and 

$G(\vec{w}) =  G_1(x_1) + \ldots + G_R(x_R)$,

\noindent then

\begin{equation} \label{RD all +separable}
\mathcal{D}(F \Vert G) \vert_W = 
\sum_{r=1}^R \mathcal{D}(F_r \Vert G_r) \vert_{W_r}.
\end{equation}

In the context of the previous example, should the batch service cost $F$ be modeled as a sum of separate costs of service $F_1, \ldots, F_R$ for each of the selections from the respective waiting lines, that $F$ becomes a completely  additive separable grading function on $W$. As such, based on formula (\ref{RD all +separable}), MRDP calls for maximizing

\begin{equation}\label{MRDP all +separable}
\mathcal{D}(F\Vert N)\vert_W = 
\sum_{r=1} ^{R} \mathcal{D}(F_r\Vert N_r)\vert_{W_R}
\end{equation}

\noindent subject to the application-specific constraints.

Furthermore, if the constraints are stated separately and independently for each of the bundled chains (waiting lines in the Queueing Theory application above) the overall maximization problems simply splits into $R$ separate independent single-chain problems. 

Otherwise, to allow for computing cost reduction, the total service cost model may sometimes be adjusted so as to reflect interdependence of costs in a convenient way. Say, one common sense way to model the batch $\vec{w}$ service cost is to present it as a sum of separate subgroups service costs and a server-specific component. That component could reflect, for example, the service batch subgroups selection, loading, post-processing, etc. Again, under the natural assumption of that cost component increasing with the batch size, the total service cost $F_0 (\vec{w})$ presents as proper grading function on $W$.
Modelling that component as linearly dependent on the batch size:

$F_0 (\vec{w}) = a+bN(\vec{w}), a>0, b>0$,

\noindent using the aforementioned linear properties of Relative Divergence, the MRDP problem simply leads to the same (\ref{MRDP all +separable}) where 
$F_r (w_r)$ is replaced by $F_r + b N_r(w_r)$.

\section{MRDP for Parameter-Controlled Grading Functions} 
\label{MRDP for controlled GFs}

Using MRDP, in some applications it has to be done by choosing values of certain controlling "intrinsic" variables ("parameters").

A representative example arises again in the context of the Queuing Theory where $R$ streams of incoming servers stop by the service station to pick up a group of customers of their designated types from the waiting line $W_r$. The server's capacity is $n_r$, and the cost of serving a group of $i_r \leq n_r$ customers by the $r$-th kind server, $F_r(i_r)$, is naturally assumed to increase with $i_r$. 

From the point of view of an observer, an assumption can be made that the incoming servers are specialized to serve customers of just one type $r = 1,\ldots,R$ of, modeled as a random variable with stationary probabilities $p_r$. The expected cost of service of a group $\vec{i} = [i_1, \ldots, i_R]$, \quad $F(\vec{i}) = \sum_{r=1}^R p_r F_r (i_r)$ is, therefore, a grading function on the chain bundle $W = \{ \vec{i} \}$.

Given the information on the service costs of groups of each type, the "most reasonable" probability distribution of the server type can be obtained by applying MRDP to $F(\vec{i})$. Using (\ref{scaling F}) and (\ref{MRDP all +separable}), it requires finding such nonnegative $p_1, \ldots,p_R, \quad p_1+ \ldots + p_R = 1, \quad $ that would maximize 

\begin{equation}\label{probs}
\sum_{r=1}^R [p_r D_r - (M_r - m_r)p_r\ln{p_r}]
\end{equation},

 \noindent where $D_r = \mathcal{D}(F_r \Vert N_r) \vert_{W_r} $.

 Using here Lagrange multipliers with respect to variables $p_1, \ldots, p_R$, it follows that at a domain-interior point of maximum, in addition to summing up to 1, they must satisfy the following system of equations:

 $D_r - (M_r - m_r) (\ln{p_r} +1) - \lambda = 0, \quad r = 1,\ldots, R$.

 Solving the system, it follows that

 $p_r = e^{[-1 + \frac{D_r - \lambda}{M_r - m_r}]}, \quad r = 1, \ldots, R$,

 \noindent with the Lagrange multiplier $\lambda$ to be determined from the condition $\quad p_1 + \ldots+ p_R = 1$.

 In some cases the service contract between the pick-up server and the queuing station may specify the same "cost" (say, the server loading time at the station) values interval) for all server types .
 
$M_r - m_r = M-m, \quad \forall r$. In that case, it follows that

 $p_r = cq^{D_r}, \quad r=1, \ldots R$,

 \noindent where $\ln{q} = (M-m)^{-1}, \quad c^{-1} = \sum_{r=1}^R q^{D_r}$.

\section{Conclusion}
\label{end}

In \cite{AOME} we introduced the concept of Relative Divergence (RD) of grading function on totally ordered sets (chains). Based on that, RD and Maximum Relative Divergence Principle (MRDP) were extended here to partially ordered sets that are direct products (bundles) of chains. 

In particular, a number of special cases arising from representative Operations Research applications were studied here. It was shown that: 

1. Similar to Maximum Entropy Principle (MEP), MRDP-based analysis leads to "most reasonable", "common sense" solutions akin to analogous Probability Theory problems;

2. A new method emerges to analyze some Operations Research applications (such as group service in Queuing Theory, resource distribution under constraints, etc.) which can be stated as MRDP problems on a bundle of chains.

3. In those applications, the computational cost of arising MRDP problems reduces strongly (and may lead to direct results) when the application of interest can be modeled in some special ways so that the grading functions emerging there can be assumed to have convenient form. In particular, it is done here when those grading functions are "height-dependent", or "additively separable", or "parameter-controlled".

\end{document}